\documentclass{article}
\usepackage{amssymb,amsfonts}

\newtheorem{theorem}{Theorem}[section]

\newcommand{\twomat}[4]{\left(\begin{array}{cc}#1&#2\\#3&#4\end{array}\right)}

\newcommand{\qed}{\hfill$\square$\par\vskip12pt}

\newcommand{\be}{\begin{equation}}
\newcommand{\ee}{\end{equation}}
\newcommand{\bea}{\begin{eqnarray}}
\newcommand{\eea}{\end{eqnarray}}
\newcommand{\beas}{\begin{eqnarray*}}
\newcommand{\eeas}{\end{eqnarray*}}

\newtheorem{proposition}{Proposition}

\newcount\minute
\newcount\hour
\def\currenttime{%
    \minute\time
    \hour\minute
    \divide\hour60
    \the\hour:\multiply\hour60\advance\minute-\hour\the\minute}
\begin{document}
\title{A Characterisation of Anti-L\"owner Functions}
\author{Koenraad M.R. Audenaert\\
Department of Mathematics,\\ Royal Holloway, University of London,\\ Egham TW20 0EX, United Kingdom}

%\subjclass[2000]{Primary 15A60}
%\commby{Marius Junge}
%\date{\today}
%\keywords{Matrix monotone functions, L\"owner matrices, Lyapunov equation}
%------------------------------------------------------------------ ABSTRACT
\maketitle
\begin{abstract}
According to a celebrated result by L\"owner, a real-valued function $f$
is operator monotone if and only
if its L\"owner matrix, which is the matrix of divided differences
$L_f=\left(\frac{f(x_i)-f(x_j)}{x_i-x_j}\right)_{i,j=1}^N$,
is positive semidefinite for every integer $N>0$ and any choice of $x_1,x_2,\ldots,x_N$.
In this paper we answer a question of R.~Bhatia, who asked for a characterisation of
real-valued functions $g$ defined on $(0,+\infty)$ for which the
matrix of divided sums $K_g=\left(\frac{g(x_i)+g(x_j)}{x_i+x_j}\right)_{i,j=1}^N$,
which we call its anti-L\"owner matrix, is positive semidefinite
for every integer $N>0$ and any choice of $x_1,x_2,\ldots,x_N\in(0,+\infty)$.
Such functions, which we call anti-L\"owner functions, have applications in the
theory of Lyapunov-type equations.
\end{abstract}

%%%%%%%%%%%%%%%%%%%%%%%%%%%%%%%%%%%%%%%%%%%%%%%%%%%%%%%%%%%%%%%%%%%%%%%%%%%%%%%%%%%%%%%%%%%%%
\section{Introduction\label{sec1}}
A real-valued function defined on an interval $(a,b)$ is called \textit{matrix monotone of order $N$} if for any pair
$A,B$ of $N\times N$ Hermitian matrices with spectrum in $(a,b)$
the implication $A\le B\Longrightarrow f(A)\le f(B)$ holds, i.e.\ $f$ preserves the positive semidefinite ordering.
A function is called \textit{operator monotone}
if it is matrix monotone of every order.

One of the central objects in the theory of matrix monotone functions is the so-called L\"owner matrix.
Given any integer $N>1$, and any set of $N$ finite, distinct real numbers $x_i$ in $(a,b)$, one constructs a
L\"owner matrix of $f$ as the $N\times N$ matrix $L_f$ of
divided differences
$$
L_f :=
\left(
\frac{f(x_i)-f(x_j)}{x_i-x_j}
\right)_{i,j=1}^N.
$$
For the diagonal elements, $i=j$, a limit has to be taken, so that
the diagonal elements are given by the first derivatives
$f'(x_i)$.
(A necessary condition for $f$ being matrix
monotone of order at least 2 is that it should be continuous,
in fact even continuously differentiable (\cite{donoghue} p.\ 79),
hence its first derivative should exist. For $N=1$, this is strictly speaking not needed, but matrix monotonicity
then reduces to ordinary monotonicity anyway.)
According to a celebrated result by L\"owner, $f$ is a matrix monotone function on $(a,b)$ of order $N$
if and only if any $N\times N$ L\"owner matrix $L_f$ is positive semidefinite,
for any choice of $x_i$ in $(a,b)$.
For a thorough introduction to matrix monotone functions we refer to the monograph \cite{donoghue},
and to \cite{bhatia} for a more concise introduction.

In (\cite{bhatia2}, p.\ 195) R.~Bhatia raised the question whether there is a good characterisation
of real-valued functions
$g(x)$ defined on $(a,b)$, with $a\ge 0$, for which every matrix of the form
$$
K_g:=\left(\frac{g(x_i)+g(x_j)}{x_i+x_j}\right)_{i,j=1}^N,
$$
is positive semidefinite, with $x_i$ distinct real numbers in $(a,b)$.
That is, $K_g$ is akin to a L\"owner matrix,
but has the minus signs replaced by plus signs.
In this paper, we'll call these matrices \textit{anti-L\"owner matrices},
and functions for which all $N\times N$ anti-L\"owner matrices are positive semidefinite
will be called \textit{anti-L\"owner functions of order $N$}.
Likewise, we call functions \textit{anti-L\"owner functions} if they satisfy this positivity
criterion for all values of $N$.

It goes without saying that to be anti-L\"owner
$g$ must first of all be non-negative, as can be seen from the trivial case $N=1$.
For $N=1$ this is already the complete answer; to avoid trivialities we will henceforth assume $N\ge 2$.
It is also straightforward to show that for $N\ge2$, $g$ must be continuous, similar to matrix monotone
functions of order $N\ge 2$; see Proposition \ref{prop:cont} below.

It has already been known for some time
that every non-negative operator monotone function on $(0,+\infty)$
is an anti-L\"owner function, and so is every non-negative operator monotone decreasing function \cite{Kwong}.
This easily follows (see, e.g.\ Theorem 1 in \cite{Kwong})
from exploiting the well-known integral representation \cite{bhatia}
\be
f(x)=\alpha+\beta x+\int_0^\infty \frac{x}{t+x}d\mu(t)\label{eq:intmono}
\ee
for non-negative operator monotone functions on $(0,+\infty)$,
with $\alpha,\beta\ge0$\footnote{Note that
for negative $\alpha$, $f$ is also operator monotone, but positivity of $f$ requires $\alpha\ge0$.}
and $\mu$ a positive Borel measure on $(0,+\infty)$ such that the given integral converges.
The statement for non-negative operator monotone decreasing functions follows easily from this by noting that
the function $f(x)$ is anti-L\"owner if and only if $1/f(x)$ is anti-L\"owner too.
%%%%%%%%%%%%%%%%%%%%%%%%%%%%%%%%%%%%%%%%%%%%%%%%%%%%%%%
\section{Main results}
In this paper, we obtain a complete answer to Bhatia's question:
\begin{theorem}\label{th:1}
Let $g$ be a real-valued function $g$, defined and finite on $(a,b)$, with $0\le a<b$.
Let $N$ be any integer, at least 2.
If $g$ is an anti-L\"owner function of order $2N$ on $(a,b)$ then
$x\mapsto g(\sqrt{x})\sqrt{x}$ is a non-negative matrix monotone function of order $N$ on $(a^2,b^2)$,
and $x\mapsto g(\sqrt{x})/\sqrt{x}$ is a non-negative matrix monotone decreasing
function of order $N$ on $(a^2,b^2)$.
\end{theorem}

Theorem \ref{th:1} has applications in the study of Lyapunov-type equations
and also answers a question by Kwong \cite{Kwong},
who studied conditions on the function $g$ such that the solution $X$ of equation $AX+XA=g(A)B+Bg(A)$ is
positive definite
for all positive definite $A$ and $B$. Kwong pointed out in \cite{Kwong2} that it suffices to consider matrices $A$
that are diagonal, with $B$ equal to the all-ones matrix ($B_{ij}=1$).
In that case the solution of the equation reduces to $X$
being an anti-L\"owner matrix with the $x_i$ equal to the diagonal elements of $A$.
Thus, our Theorem \ref{th:1} also yields the answer to Kwong's question.

As a direct corollary of Theorem \ref{th:1}, we immediately get an integral representation for anti-L\"owner
functions $g(x)$ (of all orders) on $(0,+\infty)$. From equation (\ref{eq:intmono}) we obtain
$g(\sqrt{x})=\alpha/\sqrt{x}+\beta\sqrt{x}+\int_0^\infty \frac{\sqrt{x}}{t+x}d\mu(t)$, hence
\be
g(x) = \frac{\alpha}{x}+\beta x+\int_0^\infty \frac{x}{t+x^2}d\mu(t),\label{eq:intanti}
\ee
with $\alpha,\beta\ge0$ and $\mu$ a positive Borel measure such that the integral exists.

Within this restricted setting,
the sufficiency part of Theorem \ref{th:1} is easy to prove, as it suffices to check each of the terms in
the integral representation (\ref{eq:intanti}).
To wit, one only needs to prove that the functions
$g(x)=x$ and $g(x)=x/(t+x^2)$ (for $t\ge0$) are anti-L\"owner, as all other functions
concerned are positive linear combinations of these extremal functions.
This is trivial for $g(x)=x$, because then $K_g=(1)_{i,j}$, which is clearly positive semidefinite (and rank 1).
Secondly, for $g(x)=x/(t+x^2)$, we have
\beas
K_g &=& \left(\frac{x_i/(t+x_i^2)+x_j/(t+x_j^2)}{x_i+x_j}\right)_{i,j}  \\
&=& \left(\frac{x_i(t+x_j^2) + x_j(t+x_i^2)}{(t+x_i^2)(x_i+x_j)(t+x_j^2)}\right)_{i,j} \\
&=& \left(\frac{t+x_i x_j}{(t+x_i^2)(t+x_j^2)}\right)_{i,j},
\eeas
which is congruent to the matrix $(t+x_i x_j)_{i,j}$ and therefore positive semidefinite as well
(and in general rank 2).

The hard part is to prove necessity, i.e.\ that there are no other anti-L\"owner functions than those with the
given integral representation (\ref{eq:intanti}). Furthermore, there seems to be no obvious approach even to the
sufficiency part in the more general setting of fixed $N$ where no integral representation is known.
An important observation that shows the way out, however, is hidden in the very
statement of Theorem \ref{th:1}, as it hints at a one-to-one correspondence between
anti-L\"owner functions and non-negative operator monotone functions. This is no coincidence,
and our method of proof will exploit an even deeper correspondence
between L\"owner matrices and anti-L\"owner matrices, which is made apparent in Theorem \ref{th:2} below.
This is good news, as there will be no need to develop from scratch
a completely new theory in parallel with L\"owner's.

\begin{theorem}\label{th:2}
Let $N$ be any integer and $x_1,\ldots,x_N$ a sequence of distinct positive real numbers
contained in the interval $(a,b)$, $0\le a<b$.
For any continuous real-valued function $g$ defined on $(a,b)$,
let $L$ and $K$ be its L\"owner and anti-L\"owner matrix of order $N$ on
the given points $x_1,\ldots,x_N$, respectively, and let $K_ij$ be the matrix
$$
K_{ij}=\left[\frac{g(x_k+i\epsilon)+g(x_l+j\epsilon)}{(x_k+i\epsilon)+(x_l+j\epsilon)}\right]_{k,l=1}^N
$$
and let $L_{ij}$ be the matrix
$$
L_{ij}=\left[\frac{g(x_k+i\epsilon)-g(x_l+j\epsilon)}{(x_k+i\epsilon)-(x_l+j\epsilon)}\right]_{k,l=1}^N.
$$
Then the following are equivalent:
\begin{enumerate}
\item the $2\times 2$ block matrix $\twomat{K_{00}}{K_{01}}{K_{10}}{K_{11}}$ is positive semidefinite;
\item the $2\times 2$ block matrix $\twomat{K_{00}}{L_{01}}{L_{10}}{K_{11}}$ is positive semidefinite.
\end{enumerate}
\end{theorem}

In the remainder of this paper we present the proofs of these theorems.
%%%%%%%%%%%%%%%%%%%%%%%%%%%%%%%%%%%%%%%%%%%%%%%%%%%
\section{Proofs}
We start with a simple, but nevertheless essential proposition.
\begin{proposition}\label{prop:cont}
Let $g$ be a positive real-valued function on $(a,b)$, with $0\le a<b$.
If $g$ is an anti-L\"owner function of order at least $2$, then $g$ is continuous.
\end{proposition}
\textit{Proof.}
This follows from consideration of the $2\times 2$ anti-L\"owner matrices in the
points $x_1=x,x_2=x+\epsilon$ ($0\le a<x<b$) and letting $\epsilon$ tend to 0.
Positive semidefiniteness of the anti-L\"owner matrix
requires non-negativity of its determinant:
$$g(x)g(x+\epsilon)/x(x+\epsilon)-(g(x)+g(x+\epsilon))^2/(2x+\epsilon)^2\ge 0.$$
After some calculation, one finds that this requires
$|(g(x+\epsilon)-g(x))/\epsilon|\le g(x)/x$ for all $\epsilon>0$,
whence the derivative of $g$ should exist and be bounded on any bounded
closed interval in $(a,b)$.
\qed

The main technical result on which our proof is based is the following proposition.
\begin{proposition}\label{prop:1}
Fix an integer $N$.
Let $g=(g_1,\ldots,g_N)$ and $x=(x_1,\ldots,x_N)$ be positive vectors, where all $x_i$ are distinct,
and $s=(s_1,\ldots,s_N)$ a real vector with $s_i=\pm 1$.
Then the sign of $\det Z_N$, where
$$
Z_N=\left(\frac{s_i g_i+s_j g_j}{s_i x_i+s_j x_j}\right)_{i,j=1}^N,
$$
is independent of the signs of the $s_i$'s.
\end{proposition}
As an illustration of this proposition, we will prove the easiest non-trivial case $N=2$
(the case $N=1$ is trivial as $s_1$ cancels out entirely).
The given determinant is
\beas
\det Z_2 &=& \det\twomat{\frac{g_1}{x_1}}{\frac{s_1 g_1+s_2 g_2}{s_1 x_1+s_2 x_2}}
{\frac{s_1 g_1+s_2 g_2}{s_1 x_1+s_2 x_2}}{\frac{g_2}{x_2}}\\
&=& \frac{g_1 g_2}{x_1 x_2}-\frac{(s_1 g_1+s_2 g_2)^2}{(s_1 x_1+s_2 x_2)^2} \\
&=& \frac{g_1 g_2(s_1^2 x_1^2+s_2^2 x_2^2) - x_1 x_2(s_1^2 g_1^2+s_2^2 g_2^2)}{x_1 x_2(s_1 x_1+s_2 x_2)^2}\\
&=& \frac{g_1 g_2(x_1^2+x_2^2) - x_1 x_2(g_1^2+g_2^2)}{x_1 x_2(s_1 x_1+s_2 x_2)^2}.
\eeas
One sees that the numerator is independent of the signs of the $s_i$'s,
while the denominator is always positive.
Hence, the sign of this determinant is independent of the signs of the $s_i$'s.

For small values of $N$ one can easily verify that the determinant can always be written
as a rational function where the numerator is a polynomial in which the $s_i$'s only appear to even powers,
and where the denominator is always positive. This observation provided the inspiration
for the following simple proof of Proposition \ref{prop:1} (for every value of $N$).

\textit{Proof of Proposition \ref{prop:1}.}\\
Clearly, once we prove that the sign of $\det Z_N$ does not change under a single sign change of $s_i$,
the general statement of the proposition follows, by changing the signs of the $s_i$'s one by one.
W.l.o.g.\ we consider sign changes of $s_1$ only.

The idea of the proof is to apply a \textit{partial} Gaussian elimination on $Z_N$, only bringing its first column in
upper-triangular form.
For each $i>1$ we subtract $\frac{x_1}{g_1}\,\,\frac{s_1 g_1+g_i}{s_1 x_1+x_i}$ times row 1 from row $i$.
As is well-known, this operation does not change the determinant.
The resulting matrix is of the form
$$
Z_N' =\left(
\begin{array}{cc}
\frac{g_1}{x_1} & b \\
0 & X
\end{array}
\right)
$$
where $b$ is the first row of $Z_N$ (except its element $(1,1)$)
and $X$ is an $(N-1)\times(N-1)$ matrix with elements ($i,j>1$)
\bea
X_{i,j} &=& \frac{g_i+g_j}{x_i+x_j}-\frac{x_1}{g_1}\,\,\frac{s_1 g_1+g_i}{s_1 x_1+x_i}\,\,
\frac{s_1 g_1+g_j}{s_1 x_1+x_j} \nonumber\\
&=& \frac{g_1(g_i+g_j)(x_1^2+x_i x_j)-x_1(x_i+x_j)(g_1^2+g_i g_j)}{g_1(x_i+x_j)(s_1 x_1+x_i)(s_1 x_1+x_j)}.
\label{eq:XIJ}
\eea
In the last line we have used the fact that $s_1^2=1$.

From expression (\ref{eq:XIJ})
it is clear that $X$ can be written as a matrix product $X=DYD$,
where $D$ is a diagonal matrix with diagonal elements
$1/(s_1 x_1+x_i)$ ($i>1$), and where
$Y$ is independent of $s_1$.
It follows that the determinant of $Z_N$ is given by
$$
\det Z_N = \frac{g_1}{x_1} \det(DYD) = \frac{g_1}{x_1} \det(Y) \det(D)^2.
$$
As the only factor that depends on the sign of $s_1$ appears to even power, and is therefore non-negative,
we have proven that the sign of $\det Z_N$ does not depend on the sign of $s_1$.

In a similar way, we can show that the sign of
$\det Z_N$ does not depend on any of the signs of the $s_i$'s. This ends the proof.
\qed

\textit{Proof of Theorem \ref{th:2}.}
It is a simple corollary of Proposition \ref{prop:1} that, under the conditions stated,
the positive semidefiniteness of $Z_N$
for a given choice of signs of the $s_i$'s implies positive semidefiniteness for any other choice.
Indeed, according to Sylvester's criterion, a symmetric matrix is
%positive definite if and only if its \textit{leading} principal minors are positive, and it is
positive semidefinite if and only if all its principal
minors are non-negative. In the case of $Z_N$, the principal $k\times k$ minors are determinants of the form
$\det Z_k$ ($k=1,\ldots,N$) and according to Proposition \ref{prop:1}, the signs of these determinants
are independent of the signs of the $s_i$'s appearing in them.

Consider now, in particular, the case where $N$ is even, say $N=2n$, and the $x_i$ are given by
$$(x_1,\ldots,x_n,x_{n+1},\ldots,x_{2n}) = (y_1,\ldots,y_n,y_1+\epsilon,\ldots,y_n+\epsilon)$$
for any positive $\epsilon$ small enough such that no two $x_i$ ever become equal
when $\epsilon$ tends to $0$. Let also $g_i = g(x_i)$.

We will consider two choices for the $s_i$.
Firstly, we set all $s_i=+1$. We then get the matrix
$$
K' = \twomat{K_{00}}{K_{01}}{K_{10}}{K_{11}}.
$$
Secondly, with $s_i=+1$ for $i\le n$ and $s_i=-1$ for $i>n$, we instead get
$$
K''=\twomat{K_{00}}{L_{01}}{L_{10}}{K_{11}}.
$$

As, according to Proposition \ref{prop:1}, these matrices have the same signature (same signs of the corresponding
principal minors) this yields the
equivalence of Theorem \ref{th:2}.
\qed

It is now an easy matter to prove Theorem \ref{th:1}.

\textit{Proof of Theorem \ref{th:1}.}\\
Let $N$ be a fixed integer, at least 2.
By Proposition \ref{prop:cont}, if $g$ is an anti-L\"owner function of order at least 2, then
$g$ is continuous.
Conversely, if the function $x\mapsto g(\sqrt{x}) \sqrt{x}$ is a
non-negative matrix monotone function of order $N$, then surely $g$ must be continuous too.
Thus, in any case, Theorem \ref{th:2} applies to $g$.

Let $g$ be an anti-L\"owner function of order $2N$. Thus $K'$ is positive. In the limit $\epsilon\to0$
we then find that $K_g+L_g$ and $K_g-L_g$ are positive,
due to Theorem \ref{th:2}.
A simple calculation shows that $K_g+L_g$ is equal to
\beas
K_g+L_g &=& \left(\frac{g(x_i)+g(x_j)}{x_i+x_j}+\frac{g(x_i)-g(x_j)}{x_i-x_j}\right)_{i,j=1}^N \\
&=& 2\left(\frac{x_i g(x_i)-x_j g(x_j)}{x_i^2-x_j^2}\right)_{i,j=1}^N,
\eeas
which is (up to an irrelevant factor of 2)
the L\"owner matrix of the function $x\mapsto g(\sqrt{x}) \sqrt{x}$ in the points $x_i^2$.
Hence, 
the function $x\mapsto g(\sqrt{x}) \sqrt{x}$ is a
non-negative matrix monotone function of order $N$ on $(a^2,b^2)$.

In a similar way we find
\beas
K_g-L_g &=& -2\left(x_i \,\,\frac{g(x_i)/x_i-g(x_j)/x_j}{x_i^2-x_j^2}\,\,x_j\right)_{i,j=1}^N,
\eeas
which shows that
the function $x\mapsto g(\sqrt{x})/ \sqrt{x}$ is a non-negative matrix monotone \textit{decreasing}
function of order $N$ on $(a^2,b^2)$.
\qed

%-----------------------------------------------------------------------
\section*{Acknowledgments}
The author gratefully acknowledges the hospitality of the University of Ulm, Germany, and
of the Institut Mittag-Leffler, Stockholm,
where parts of this work have been done.
Also many thanks to R.~Bhatia and F.~Hiai, for comments on an earlier version
of the manuscript.
%------------------------------------------------------------- BIBLIOGRAPHY
\bibliographystyle{amsplain}

%%%%%%%%%%%%%%%%%%%%%%%%%%%%%%%%%%%%%%%%%%%%%%%%%%%%%%%%%%%%%%%%%%%
\end{document}